\documentclass[12pt]{article}

\usepackage{epsfig}
\usepackage{amsmath,ams}
\usepackage{soul,color}
\usepackage{bm}
\usepackage{epsfig}
\usepackage{natbib}
\newcommand{\tr}{^{\prime}}
\def\b#1{\mbox{\boldmath $#1$}}    

\newcommand{\ot}{\mbox{$\:\otimes \:$}}
\renewcommand{\th}{\theta}

\newcommand{\la}{\lambda}
\newcommand{\si}{\sigma}

\newtheorem{theo}{Theorem}
\newtheorem{example}{Example}

\usepackage{bm}
\usepackage{epsfig}

\begin{document}
\title{\vspace*{-2cm}An alternative to the Baum-Welch recursions\\ for hidden Markov models}
\author{Francesco Bartolucci\footnote{Department of Economics, Finance and
Statistics, University of Perugia, Via A. Pascoli, 20, 06123
Perugia. E-mail: bart@stat.unipg.it} \footnote{I acknowledge the financial support
from the ``Einaudi for Economics and Finance'' (Rome - IT)}}
\date{}
\maketitle
\vspace*{-0.5cm}
\begin{abstract}
We develop a recursion for hidden Markov model of any order $h$, which allows us to obtain
the posterior distribution of the latent state at every occasion, given the previous
$h$ states and the observed data. With respect to the well-known Baum-Welch recursions,
the proposed recursion has the advantage of being more direct to use and, in particular,
of not requiring dummy renormalizations to avoid numerical problems. We also
show how this recursion may be expressed in matrix notation, so as to allow for
an efficient implementation, and how it may be used to obtain the manifest distribution
of the observed data and for parameter estimation
within the Expectation-Maximization algorithm. 
The approach is illustrated by an application to financial data
which is focused on the study of the dynamics of the volatility level of log-returns.
\end{abstract}\newpage
\section{Introduction}
Hidden Markov (HM) models have become a popular statistical tool for the analysis of
data having a time-series structure; for an up-to-date review see \cite{zucchini:macdonald:2009}.
These models have also found great interest for the analysis of longitudinal data,
where independent short time series are observed for typically many
statistical units; for a review
see \cite{bart:farc:pen:10}. HM models are based on the assumption that the 
observable random variables, corresponding to the different time occasions, are conditionally
independent given an unobservable (or latent) process, which follows a Markov chain.
Usually, this Markov chain is assumed to be of first order and time homogenous, so that
the transition probabilities are time invariant.

A fundamental tool of inference for HM models is represented by forward-backward recursions
of Baum and Welch \citep[see][]{baum:et:al:1970,welch:2003}. For a first-order HM model,
these recursions allow us to compute the manifest probability 
(or density) of the observed sequence of data and to obtain the posterior distribution of 
every latent state and of every pair of consecutive latent states given these data. 
Through this recursion is then possible to implement an Expectation-Maximization (EM) algorithm
\citep{baum:et:al:1970,demp:lair:rubi:77} for maximum likelihood estimation of the parameters
and to perform local decoding \citep{juan:rabi:91}, that is to find the most likely state
at every occasion, given the observed data.
Despite its popularity, the Baum-Welch recursions
may suffer from numerical problems due to the fact
that certain probabilities may become negligible. This problem typically requires to implement
dummy renormalizations; see \cite{Scott:02} for further comments and \cite{lyst:hugh:02}
for an alternative solution in dealing with the manifest distribution of the observed
data.

In a rather recent paper, \cite{bart:besag:2002} proposed a probabilistic
result to obtain the marginal distribution of a random variable in
Markov random field model and mentioned that this result may be also used
for HM models, providing an example for a first-order and a second-order HM model.
Developing the intuition of \cite{bart:besag:2002}, 
in this paper we propose a general recursion
to deal with HM models of any order $h$. This recursion allows us to
obtain the posterior distribution of every latent state given the previous
$h$ states and the observed data. With respect to the Baum-Welch recursions,
the proposed recursion has the advantage of being more direct to use,
especially with higher-order HM models. 
Moreover, it does not require dummy renormalizations. 

We show how the proposed recursion may be used to obtain the manifest distribution of
the observed data and the required posterior probabilities to implement the EM algorithm
for parameter estimation. Moreover, the recursion may be directly used for local decoding 
and for prediction. In order to allow for an efficient implementation, 
we also express the proposed result in matrix notation. Such an implementation in the
{\sc R} language is available to the reader upon request.

The remainder of the paper is organized as follows. In the following section we briefly
review HM models and the Baum-Welch recursion. The proposed recursion is illustrated
in Section 3, whereas in Section 4 we illustrate its use for maximum likelihood
estimation, local decoding, and prediction. Finally,
in Section 5 we provide an illustration by an application based on an HM version
of the stochastic volatility (SV) model for financial data \citep{taylor:2005},
in which we assume the existence of discrete levels of volatility.
\section{Preliminaries}
Consider a sequence of $T$ {\em manifest random variables} $Y_1,\ldots,Y_T$ which
are collected in the vector $\b Y$. A hidden Markov
(HM) model assumes that these random variables are conditional independent given
the {\em unobservable random variables} $U_1,\ldots,U_T$ which follow a Markov chain
with $k$ states.
We consider in particular a Markov chain of order $h$ so that
\[
p(u_t|u_1,\ldots,u_{t-1})=
p(u_t|u_{t-h},\ldots,u_{t-1}),\quad t=h+1,\ldots,T,
\]
where we use the notation $p(u_t|u_1,\ldots,u_{t-1}) = P(U_t=u_t|U_1=u_1,\ldots,U_{t-1}=u_{t-1})$.
A similar notation will be adopted throughout the paper
to denote probability mass functions, in a way that will be clear from the context.
It is also assumed that every $Y_t$ depends on the latent process only through
$U_t$ and then by $f(y_t|u_t)$ we denote the probability mass (or density) function
of this distribution.

The specific HM model adopted in an application is based on assumptions
on the above transition probabilities, such as that these probabilities are
time homogeneous. These assumptions may also concern the distribution of each
response variable given the corresponding latent variable. The specific formulation
may also involve covariates, if available. In this section, however, we remain in 
the general context described above and base most results on the unspecified transition
probability function $p(u_t|u_{\max(t-h,1)},\ldots,u_{t-1})$ and the conditional
response probability (or density) function $f(y_t|u_t)$. Note that, in 
denoting the transition probabilities, we use the index $\max(t-h,1)$ in order
to have a notation that is suitable even for $t<h$. Obviously, when $t=1$, the
conditioning argument in these probabilities vanishes and they reduce to 
initial probabilities of type $p(u_1)$.

The following example clarifies a possible formulation of an HM model for time-series data.
For other examples in the context of longitudinal data see \cite{bart:farc:pen:10}.

\begin{example}\label{ex1}
Consider an HM version of the SV model 
for financial data \citep{taylor:2005},
which is based on the assumption that, given $U_t$, the log-return
$Y_t$ has a normal distribution with mean 0 and variance depending on $U_t$.
In particular, we assume that
\[
f(y_t|u_t)=\frac{1}{\sqrt{2\pi\si^2_{u_t}}}\exp\left[
-\frac{1}{2}\left(\frac{y_t}{\si_{u_t}}\right)^2
\right],
\]
where $\si_v$, $v=1,\ldots,k$, are volatility levels associated to the different
latent states.
We also assume that the underlying Markov chain is of order $h$ and is time-homogenous, 
so that, for all $t>h$, we have
\[
p(u_t|u_{t-h},\ldots,u_{t-1})=\pi_{u_{t-h},\ldots,u_t},
\] 
where $\pi_{v_1,\ldots,v_{h+1}}$, $v_1,\ldots,v_{h+1}=1,\ldots,k$, are common 
transition probabilities to be estimated together with
$\si_1,\ldots,\si_k$. Other parameters to be estimated are the initial
and transition probabilities for $t\leq h$. These parameters are denoted
by
\[
\la_{t,u_{\max(t-h,1)},\ldots,u_t}=p(u_t|u_{\max(t-h,1)},\ldots,u_{t-1}).
\]
Overall, taking into account that the initial probabilities are such that
$\sum_{u_1}\la_{1,u_1}=1$ and similar constraints hold for all transition
probabilities, the number of free parameters is
\begin{equation}
\#par = \underbrace{k}_{\si_v}+
\underbrace{(k-1)\sum_{t=1}^{h-1} k^{t-1}}_{\la_{t,u_{\max(t-h,1)},\ldots,u_t}}
+ \underbrace{(k-1)k^h}_{\pi_{v_1,\ldots,v_{h+1}}}.\label{eq:npar}
\end{equation}
\end{example}\vspace*{0.5cm}

It has to be clear that the same modeling framework described above may be adopted with
longitudinal data in which we observe short sequences of data for $n$ sample
units, which are usually assumed to be independent. However, we do not
explicitly consider the case of longitudinal data
since the theory that will be developed easily apply
to this case as well.

In order to efficiently compute the probability (or the density) of an
observed sequence of $T$ observations, collected in the vector $\b y=(y_1,\ldots,y_T)$,
Baum and Welch \citep{baum:et:al:1970,welch:2003} proposed the following
forward recursion for a first-order HM model:
\begin{equation}
f(u_t,\b y_{\leq t})=\sum_{u_{t-1}}f(u_{t-1},\b y_{\leq t-1})
p(u_t|u_{t-1})f(y_t|u_t),\quad t=2,\ldots,T,\label{eq:forward_recusion}
\end{equation}
where $\b y_{\leq t}=(y_1,\ldots,y_t)$. This recursion is
initialized with $f(u_1,y_1)=p(u_1)f(y_1|u_1)$ and, in the end,
we obtain the manifest probability (or density) function of $\b y$ as
\[
f(\b y)=\sum_{u_t}f(u_t,\b y).
\]

Moreover, Baum and Welch introduced the backward recursion
\begin{equation}
f(\b y_{>t}|u_t)=\sum_{u_{t+1}}f(\b y_{>t+1}|u_{t+1})p(u_{t+1}|u_t)f(y_{t+1}|u_{t+1}),
\quad t=1,\ldots,T-1,
\label{eq:backward_recursion}
\end{equation}
where $\b y_{> t}=(y_{t+1},\ldots,y_T)$, which is initialized with $f(\b y_{>T}|u_t)=1$.
Using this recursion, we can obtain the posterior probability of every
latent state given the observed data, that is $q(u_t|\b y)=P(U_t=u_t|\b Y=\b y)$.
In particular, we have
\[
q(u_t|\b y)=\frac{f(u_t,\b y_{\leq t})f(\b y_{>t}|u_t)}{f(\b y)},\quad t=1,\ldots,T,
\]
whereas for the posterior probability of every pair of consecutive states we have
the posterior probability
\[
q(u_{t-1},u_t|\b y)=\frac{f(u_{t-1},\b y_{\leq t})p(u_t|u_{t-1})f(y_t|u_t)
f(\b y_{>t}|u_t)}{f(\b y)},\quad t=2,\ldots,T.
\]

As mentioned above, the Baum-Welch recursions suffer from the problem of numerical 
instability due to the fact that, as $t$ increases, the probability in 
(\ref{eq:forward_recusion}) becomes negligible. The problem is evident when
$T$ is large and also affects the probabilities in (\ref{eq:backward_recursion}).
This problem requires suitable renormalizations; 
see \cite{Scott:02} for a more detailed description. 
\section{Proposed recursion}
Developing a result due to \cite{bart:besag:2002} for Markov random fields,
in this section we propose how to compute the posterior probabilities
\begin{equation}
q(u_t|u_{\max(t-h,1)},\ldots,u_{t-1},\b y),\quad t=1,\ldots,T,
\label{eq:def:p_t}
\end{equation}
that is the conditional probability of a certain realization of $U_t$,
given $U_{\max(t-h,1)},\ldots,U_{t-1}$ and a certain configuration of
responses collected in the vector $\b y$.

For last time occasion, that is when $t=T$, the above probability may be
simply computed as
\begin{equation}
q(u_T|u_{\max(T-h,1)},\ldots,u_{T-1},\b y)=
\frac{f(y_T|u_T)p(u_T|u_{\max(T-h,1)},\ldots,u_{T-1})}
{c(u_{\max(T-h,1)},\ldots,u_{T-1},y_T)},\label{eq:p_T}
\end{equation}
where $c(u_{\max(t-h,1)},\ldots,u_{T-1},y_T)$ is the normalizing
constant equal to the sum of the numerator of (\ref{eq:p_T}) 
for all the possible values of $U_T$. 

Now consider the following Theorem that allows us to compute the
conditional probability in (\ref{eq:def:p_t}) for $t$ smaller than $T$
and is related to Theorem 1 of \cite{bart:besag:2002}.

\begin{theo}
We have that
\begin{eqnarray}
&&q(u_t|u_{\max(t-h,1)},\ldots,u_{t-1},u_{t+1},\ldots,u_{t+j},\b y)=\nonumber\\
&&=
\left[
\sum_{u_{t+j+1}}
\frac{q(u_{t+j+1}|u_{\max(t+j+1-h,1)},\ldots,u_{t+j},\b y)}
{q(u_{t}|u_{\max(t-h,1)},\ldots,u_{t-1},u_{t+1},\ldots,u_{t+j+1},\b y)}
\right]^{-1},\label{eq:theo}
\end{eqnarray}
for $t=1,\ldots,T-1$ and $j=0,\ldots,\min(h,T-t)-1$ and where the conditioning variables 
$u_{t+1},\ldots,u_{t+j}$ at lhs vanishes for $j=0$.
\end{theo}

\noindent{\bf Proof} First of all consider that the assumption
that the latent Markov process is of order $h$ implies that
\[
q(u_{t+j+1}|u_{\max(t+j+1-h,1)},\ldots,u_{t+j},\b y)=
q(u_{t+j+1}|u_{\max(t-h,1)},\ldots,u_{t+j},\b y)
\]
and then we have
\[
\frac{q(u_{t+j+1}|u_{\max(t+j+1-h,1)},\ldots,u_{t+j},\b y)}
{q(u_{t}|u_{\max(t-h,1)},\ldots,u_{t-1},u_{t+1},\ldots,u_{t+j+1},\b y)}=
\frac{p(u_{\max(t-h,1)},\ldots,u_{t-1},u_{t+1},\ldots,u_{t+j+1},\b y)}
{p(u_{\max(t-h,1)},\ldots,u_{t+j},\b y)}. 
\]
Consequently, the sum in (\ref{eq:theo}) is equal to
\[
\frac{p(u_{\max(t-h,1)},\ldots,u_{t-1},u_{t+1},\ldots,u_{t+j},\b y)}
{p(u_{\max(t-h,1)},\ldots,u_{t+j},\b y)}
\]
and the Theorem holds. $\Box$\vspace*{1cm}

On the basis of the above result, we implement a backward recursion finalized to
computing the probabilities in (\ref{eq:def:p_t}). As already mentioned, for
$t=T$ these probabilities may be directly obtained from (\ref{eq:p_T}). Then,
in reverse order for $t=1,\ldots,T-1$ we first compute the posterior probabilities
\[
q(u_{t}|u_{\max(t-h,1)},\ldots,u_{t-1},u_{t+1},\ldots,u_{t+j},\b y),
\]
with $j=\min(T-t,h)$. 
Since $U_t$ is conditionally independent of $Y_1,\ldots,Y_{t-1},Y_{t+1},\ldots,Y_T$
given $U_{\max(t-h,1)},\ldots,U_{t-1}$, $U_{t+1},\ldots,U_{t+j}$, and $Y_t$, we have that
the above probability is equal to
\begin{equation}
\frac{f(y_t|u_t)\prod_{l=0}^{j}p(u_{t+l}|u_{\max(t+l-h,1)},\ldots,u_{t+l-1})}
{c(u_{\max(t-h,1)},\ldots,u_{t-1},u_{t+1},\ldots,u_{t+j},y_t)}.
\label{eq:prob_cond2}
\end{equation}
The normalizing constant at the denominator is obtained by summing
the numerator for all possible values of $u_t$. Then we apply result (\ref{eq:theo})
from $j=\min(T-t,h)-1$ to $j=0$, so as to
recursively remove the dependence of $U_t$ on $U_{t+j+1}$ and obtaining
the target posterior probabilities $q(u_t|u_{\max(t-h,1)},\ldots,u_{t-1},\b y)$.

In order to clarify the above algorithm, we explicit 
consider below the case of a first-order 
and a second-order HM model.

\begin{example}
For a first-order model ($h=1$),
the algorithm consists of first computing the probabilities
\[
q(u_T|u_{T-1},y_T)=\frac{f(y_T|u_T)p(u_T|u_{T-1})}
{c(u_{T-1},y_T)}.
\]
Then, we for $t=1,\ldots,T-1$ we apply the rule in (\ref{eq:theo}) in reverse
order. In particular, for $T\geq 3$, we have
\[
q(u_t|u_{t-1},\b y)=
\left[
\sum_{u_{t+1}}
\frac{q(u_{t+1}|u_t,\b y)}{q(u_{t}|u_{t-1},u_{t+1},\b y)}
\right]^{-1},\quad t=2,\ldots,T-1,
\]
and
\[
q(u_1|\b y)=
\left[
\sum_{u_2}
\frac{q(u_2|u_1,\b y)}{q(u_1|u_2,\b y)}\right]^{-1},
\]
where
\begin{eqnarray*}
q(u_1|u_2,y_1)&=&\frac{f(y_1|u_1)p(u_2|u_1)}{c(u_2,y_1)},\\
q(u_t|u_{t-1},u_{t+1},y_t)&=&\frac{f(y_t|u_t)p(u_t|u_{t-1})p(u_{t+1}|u_t)}
{c(u_{t-1},u_{t+1},y_t)},\quad t=2,\ldots,T-1.
\end{eqnarray*}
\end{example}

\begin{example}
For a second-order model ($h=2$), the algorithm consists of first computing the probabilities
\[
q(u_T|u_{T-2},u_{T-1},y_T)=\frac{f(y_T|u_T)p(u_T|u_{T-2},u_{T-1})}
{c(u_{T-2},u_{T-1},y_T)}.
\]
Then, for $t=1,\ldots,T-1$ we apply the rule in
(\ref{eq:theo}) for $j=2$ (provided that $t\leq T-2$) and then for $j=1$. 
In particular, assuming that $T\geq 4$, we first compute
\[
q(u_{T-1}|u_{T-3},u_{T-2},u_T,y_T)=\frac{f(y_{T-1}|u_{T-1})p(u_{T-1}|u_{T-3},u_{T-2})
p(u_T|u_{T-2},u_{T-1})}{c(u_{T-3},u_{T-2},u_T,y_T)}
\]
and consequently
\[
q(u_{T-1}|u_{T-3},u_{T-2},\b y)=
\left[
\sum_{u_T}
\frac{q(u_T|u_{T-2},u_{T-1},\b y)}
{q(u_{T-1}|u_{T-3},u_{t-2},u_T,\b y)}
\right]^{-1}.
\]
Then in reverse order for $T=3,\ldots,T-2$, we first compute
\begin{eqnarray*}
q(u_t|u_{t-2},u_{t-1},u_{t+1},u_{t+2},y_t)&=&\frac{f(y_t|u_t)p(u_t|u_{t-2},u_{t-1})
p(u_{t+1}|u_{t-1},u_t)}{c(u_{t-2},u_{t-1},u_{t+1},u_{t+2},y_t)}\times\\
&&\hspace*{1cm}\times p(u_{t+1}|u_{t-1},u_t)p(u_{t+2}|u_t,u_{t+1}),
\end{eqnarray*}
we remove the dependence of $U_t$ on $U_{t+2}$ by computing
\[
q(u_t|u_{t-2},u_{t-1},u_{t+1},\b y)=
\left[
\sum_{u_{t+2}}
\frac{q(u_{t+2}|u_t,u_{t+1},\b y)}
{q(u_{t}|u_{t-2},u_{t-1},u_{t+1},u_{t+2},\b y)}
\right]^{-1},
\]
and finally we remove the dependence on $U_{t+1}$ by computing
\[
q(u_t|u_{t-2},u_{t-1},\b y)=
\left[
\sum_{u_{t+1}}
\frac{q(u_{t+1}|u_{t-1},u_t,\b y)}
{q(u_{t}|u_{t-2},u_{t-1},u_{t+1},\b y)}
\right]^{-1}.
\]
In the end, we use similar rules to obtain $q(u_2|u_1,\b y)$ and consequently $q(u_1|\b y)$
on the basis of
\begin{eqnarray*}
q(u_1|u_2,u_3,y_1)&=&\frac{f(y_1|u_1)p(u_2|u_1)p(u_3|u_1,u_2)}{c(u_2,u_3,y_1)},\\
q(u_2|u_1,u_3,u_4,y_2)&=&\frac{f(y_2|u_2)p(u_2|u_1)p(u_3|u_1,u_2)p(u_4|u_2,u_3)}
{c(u_2,u_3,y_2)}.
\end{eqnarray*}\vspace*{5mm}
\end{example}

A crucial point is applying these recursions is the efficient implementation. At this
regard, it is worth noting that for the first-order HM model we can express the
recursion in matrix notation and then efficiently implement it in languages such as
{\sc Matlab} and {\sc R}. Details on this are provided in Appendix.
\section{Maximum likelihood estimation using the proposed recursion}
Given a sequence of observations $y_1,\ldots,y_T$ collected in $\b y$,
the {\em model log-likelihood} is 
\begin{equation}
\ell(\b\th)=\log p(\b y)
\label{eq:log_lik}
\end{equation}
where $\b\th$ is vector collecting all model parameters. 
The structure of $\b\th$ depends on the specific
parametrization which is adopted for the conditional response distribution 
$f(y_t|u_t)$ and the initial and transition probabilities 
$p(u_t|u_{\max(t-h,1)},\ldots,u_{t-1})$. For instance, for the
HM-SV model illustrated in Example \ref{ex1}, $\b\th$ includes the initial
and transition probabilities $\la_{t,u_{\max(t-h,1)},\ldots,u_t}$ and 
$\pi_{v_1,\ldots,v_{h+1}}$ and the standard deviations $\si_v$. We recall
that, in this case, the probabilities $\pi_{v_1,\ldots,v_{h+1}}$ are common
to all $t>h$, being the underlying Markov chain time homogenous.

In the following, we show how
to compute the log-likelihood in (\ref{eq:log_lik})
and implement its maximization by the recursion developed
in the previous section. It has to be clear that the same algorithm may be used in
with longitudinal data, even in the presence of individual covariates.

First of all, for any sequence of latent states $u_1,\ldots,u_T$ collected
in $\b u$, we simply have that
\[
p(\b y) = \frac{f(\b u,\b y)}{q(\b u|\b y)}=
\frac{\prod_t f(y_t|u_t)p(u_t|u_{\max(t-h,1)},\ldots,u_{t-1})}
{\prod_t q(u_t|u_{\max(t-h,1)},\ldots,u_{t-1},\b y)},
\]
where $f(\b u,\b y)$ refers to the joint distribution 
of $U_1,\ldots,U_T$ and $Y_1,\ldots,Y_T$ and $q(\b u|\b y)$ to the posterior
distribution of $U_1,\ldots,U_T$ given $Y_1,\ldots,Y_T$.
Consequently, given an arbitrary sequence $\b u$, say that with all states
equal to 1, we compute the model log-likelihood as
\[
\ell(\b\th)=\sum_t \log
\frac{f(y_t|u_t)p(u_t|u_{\max(t-h,1)},\ldots,u_{t-1})}
{q(u_t|u_{\max(t-h,1)},\ldots,u_{t-1},\b y)}
\]
on the basis of the proposed recursion. Note that, in this way, we do not
need to use any renormalization, which are instead necessary in the Baum and Welch recursions;
see also \cite{lyst:hugh:02}.

In order to maximize $\ell(\b\th)$, we can use an Expectation-Maximization algorithm
that follows the same principle as that illustrated by \cite{baum:et:al:1970}. In
particular, this algorithm is based on the {\em complete data log-likelihood}
\begin{eqnarray}
\ell^*(\b\th) &=& \sum_t\sum_{u_t} w_{t,u_t}\log f(y_t|u_t) +\nonumber \\
&+&\sum_t\sum_{u_{\max(t-1,h)}}\cdots\sum_{u_t}
z_{t,u_{\max(t-h,1)},\ldots,u_t}\log p(u_t|u_{\max(t-h,1)},\ldots,u_{t-1}),\label{eq:comp_lk}
\end{eqnarray}
where $w_{t,u_t}$ is a dummy variable equal to 1 if the latent state at occasion
$t$ is $u_t$ and to 0 otherwise and $z_{t,u_{\max(t-h,1)},\ldots,u_t}$
is a corresponding dummy variable for the
sequence of latent states $u_{\max(t-h,1)},\ldots,u_t$, which may be expressed
through the product 
\[
z_{t,u_{\max(t-h,1)},\ldots,u_t}=w_{\max(t-h,1),u_{\max(t-h,1)}}\cdots
w_{t,u_t}.
\]

At the E-step of the EM algorithm, we need to compute the posterior expected values
of the above dummy variables given the observed data and the current value of the
parameters. In particular, we have that
\begin{eqnarray*}
\hat{w}_{t,u_t}&=&E(w_{t,u_t}|\b y)=q(u_t|\b y),\\
\hat{z}_{t,u_{\max(t-h,1)},\ldots,u_t}&=&
E(z_{t,u_{\max(t-h,1)},\ldots,u_t}|\b y)=q(u_{\max(t-h,1)},\ldots,u_t|\b y).
\end{eqnarray*}
In particular, from the proposed recursion, we directly obtain $q(u_1|\b y)$. Then,
for $t>1$, we exploit a trivial forward recursion:
\begin{eqnarray}
&&q(u_{\max(t-h,1)},\ldots,u_t|\b y)=\nonumber\\
&&=\left\{\begin{array}{ll}
q(u_t|u_{\max(t-h,1)},\ldots,u_{t-1},\b y)q(u_{\max(t-h,1)},\ldots,u_{t-1}|\b y), & t=2,\ldots,h+1,\\
q(u_t|u_{t-h},\ldots,u_{t-1},\b y)\sum_{u_{t-h-1}}q(u_{t-h-1},\ldots,u_{t-1}|\b y), & t=h+2,\ldots,T,
\end{array}
\right.\label{eq:forward_recursion2}
\end{eqnarray}
to be performed for $t=1,\ldots,T$. Then, $q(u_t|\b y)$ is computed by a 
suitable marginalization. How to formulate the above forward recursion in 
matrix notation, so as to efficiently implement it, is illustrated in Appendix.

As usual, the M-step of the EM algorithm consists of maximizing $\ell^*(\b\th)$, once
the dummy variables in (\ref{eq:comp_lk})
are substituted by the corresponding expected values obtained as above.
The following example clarify how to implement this step for a specific model.

\begin{example}
For the HM-SV model illustrated in Example \ref{ex1}, the parameters $\si_v$ are
updated at the M-step as follows:
\[
\si_v = \sqrt{\frac{\sum_t \hat{w}_{t,v}y_t^2}
{\sum_t \hat{w}_{t,v}}},\quad v=1,\ldots,k.
\]
Moreover, for the initial probabilities we have
\[
\la_{1,u_1}=\hat{w}_{1,u_1},\quad u_1=1,\ldots,k,
\]
and for the transition probabilities, we have
\[
\la_{t,u_{\max(t-h,1)},\ldots,u_t}=\frac{\hat{z}_{t,u_{\max(t-h,1)},\ldots,u_t}}
{\sum_v \hat{z}_{t,u_{\max(t-h,1)},\ldots,u_{t-1},v}},
\quad u_{\max(t-h,1)},\ldots,u_t=1,\ldots,k,
\]
for $t=1,\ldots,h$ and
\[
\pi_{v_1,\ldots,v_{h+1}}=\frac{\sum_{t>h}\hat{z}_{t,v_1,\ldots,v_{h+1}}}
{\sum_v\sum_{t>h}\hat{z}_{t,v,\ldots,v_h,v}},
\quad v_1,\ldots,v_{h+1}=1,\ldots,k.
\]
\end{example}\vspace*{5mm}

Clearly, the posterior probability obtained by the proposed recursion may also be used
for local decoding \citep{juan:rabi:91}, that is to find the most likely value $\hat{u}_t$
of the latent state $U_t$, given the observed data. In particular, $\hat{u}_t$ is
found as the value that maximize the posterior probability $q(u_t|\b y)$.

Finally, on the basis of a sequence of $h$ latent states of the type $u_{T-h+1},\ldots,u_T$,
which may be even fixed by the local decoding method,
it is possible to predict the latent state at occasion $T+1$, denoted by $\hat{u}_{T+1}$,
as the value which maximizes $q(u_{T+1}|u_{T-h+1},\ldots,u_T,\b y)$; 
this posterior probability directly
derives from the proposed recursion. We can also predict the manifest distribution
of $Y_{T+1}$ through the following finite mixture
\[
\sum_{u_{T+1}}f(y_{T+1}|u_{T+1})q(u_{T+1}|u_{T-h+1},\ldots,u_T,\b y).
\]
\section{An application}
In order to illustrate the proposed approach, we fitted the HM version of the stochastic
volatility model described in Example 1 to the SP500 data for the period from the
beginning of 2008 to the end 2011. The observed outcome is the percentage log-return
with respect to the previous closing day, so that we have $T=1007$ observations.

For the above data, we estimated the model at issue for different values of $h$
(order of the latent Markov chain) and different values of the number of $k$
(number of latent sates), by the EM algorithm outlined in the previous section. The aim
of this preliminary analysis is to check if the assumption that the volatility level
follows a first-order process
is plausible. This means that the level of volatility in a given day only depends on
that of the previous day. This hypothesis may be compared with that of a higher-order dependence,
in which the level of volatility in a given day also depends on the volatility of, say, the
previous two days.

The results of the preliminary fitting are reported in Table \ref{table1} in terms of
log-likelihood, number of parameters, computed as in (\ref{eq:npar}),
and Bayesian Information Criterion \citep[BIC; ][]{sch:78}.
Note that we also include results for the model with $h=0$, which assumes independence
between the volatility levels corresponding to different time occasions.

\begin{table}[ht]
\begin{center}
\begin{tabular}{lrrrrr}
  \hline
  & & \multicolumn4c{$k$}\\\cline{3-6}
  & $h$ & \multicolumn1c{1} & \multicolumn1c{2} & \multicolumn1c{3} & \multicolumn1c{4} \\ 
  \hline
  & 0 & -2026.60 & -1898.73 & -1887.46 & -1885.57 \\ 
log-lik.  & 1 & -2026.60 & -1819.45 & -1778.00 & -1764.06 \\ 
  & 2 & -2026.60 & -1807.69 & -1768.97 & -1746.45 \\ \hline
 &   0 & 1 & 3 & 5 & 7 \\ 
\#par      & 1 & 1 & 5 & 11 & 19 \\ 
      & 2 & 1 & 9 & 29 & 67 \\ \hline
   & 0 & 4060.12 & 3818.19 & 3809.50 & 3819.54 \\ 
BIC      & 1 & 4060.12 & 3673.48 & 3632.05 & 3659.49 \\ 
      & 2 & 4060.12 & 3677.61 & 3738.46 & 3956.18 \\ 
   \hline
\end{tabular}
\end{center}
\caption{\em Results from the preliminary fitting, in terms of maximum log-likelihood,
number of parameters, and BIC, of the HM-SV model for different values of
$h$ (latent Markov chain order) and $k$ (number of latent states).}\label{table1}
\vspace*{5mm}
\end{table}

According to BIC, the observed data supports the hypothesis of a first-order dependence
of the stochastic volatility. In fact, the smallest value of the BIC index, among
those in Table \ref{table1}, is observed for $h=1$ and $k=3$. For this model, we report
in Table \ref{table2} the estimates of the parameters of main interest.

\begin{table}[ht]
\begin{center}
\begin{tabular}{rrrrrrr}
  \hline
     &                               &&       & \multicolumn3c{$\hat{\pi}_{v_1,v_2}$} \\\cline{5-7}
 $v$ & \multicolumn1c{$\hat{\si}_v$} && $v_1$ & 
 \multicolumn1c{$v_2=1$} & \multicolumn1c{$v_2=2$} & \multicolumn1c{$v_2=3$} \\ \hline
  1 & 0.865 && 1 & 0.988 & 0.010 & 0.002 \\ 
  2 & 1.609 && 2 & 0.013 & 0.981 & 0.006 \\ 
  3 & 3.770 && 3 & 0.000 & 0.025 & 0.975 \\ 
   \hline
\end{tabular}
\end{center}
\caption{\em Estimates of the parameters $\si_v$ (levels of volatility) and $\pi_{v_1,v_2}$
(transition probabilities) under the HM-SV model with $h=1$ and $k=3$.}\label{table2}\vspace*{5mm}
\end{table}
\newpage

We then observe three distinct levels of stochastic volatility and very high persistence
in the volatility level, since the probabilities in the transition matrix in Table
\ref{table2} are very close to 1. As a comparison, we report in Table \ref{table3} the
corresponding parameter estimates under the model with $h=2$ and $k=3$.

\begin{table}[ht]
\begin{center}
\begin{tabular}{rrrrrrrr}
  \hline
     &&&&& \multicolumn3c{$\hat{\pi}_{v_1,v_2,v_3}$}\\\cline{6-8}
 $v$ & \multicolumn1c{$\hat{\si}_v$}  && $v_1$ & $v_2$ &
 \multicolumn1c{$v_3=1$} & \multicolumn1c{$v_3=2$} & \multicolumn1c{$v_3=3$} \\ 
  \hline
  1 & 0.842 && 1 & 1 & 0.979 & 0.021 & 0.000 \\ 
  2 & 1.725 && 1 & 2 & 0.909 & 0.091 & 0.000 \\ 
  3 & 4.047 && 1 & 3 & 0.585 & 0.000 & 0.415 \\ 
    &       && 2 & 1 & 0.113 & 0.873 & 0.014 \\ 
    &       && 2 & 2 & 0.027 & 0.966 & 0.007 \\ 
    &       && 2 & 3 & 1.000 & 0.000 & 0.000 \\ 
    &       && 3 & 1 & 0.000 & 0.000 & 1.000 \\ 
    &       && 3 & 2 & 0.000 & 1.000 & 0.000 \\ 
    &       && 3 & 3 & 0.000 & 0.035 & 0.965 \\ 
   \hline
\end{tabular}
\end{center}
\caption{\em Estimates of the parameters $\si_v$ (levels of volatility) and $\pi_{v_1,v_2,v_3}$
(transition probabilities) under the HM-SV model with $h=2$ and $k=3$.}
\label{table3}\vspace*{5mm}
\end{table}

We observe that the estimated levels of volatility under the second-order model
are similar to those under the first-order model. Moreover, we again note a high
persistence, in the sense that $\hat{\pi}_{v_1,v_2,v_3}$ is very close to 1 whenever
$v_1=v_2=v_3$. The estimates of these transition probabilities for $v_1\neq v_2$
seem to be less reasonable, especially when $v_1=3$ and $v_2=1$. However,
we have to consider that a jump from state 3 to state 1 is very rare and then
there is no support from the data to estimate a transition probability given
this pair of states. This confirms that the first-order model is preferable for
the data at hand and may provide more reliable estimates. In any case, the possibility
to estimate a higher order HM model, which is allowed by the proposed recursion,
is important in order to have a counterpart against which comparing the more common
first-order model.
\section*{Appendix: the recursion in matrix notation}
First of all let $\b p_t$ be the column vector of prior probabilities 
$p(u_t|u_{\max(t-h,1)},\ldots,u_{t-1})$ arranged in lexicographical order so that,
for instance, with $h=2$ and $k=2$ we have
\[
\b p_t=\left(\begin{matrix}
p(u_t=1|u_{t-2}=1,u_{t-1}=1)\cr
p(u_t=2|u_{t-2}=1,u_{t-1}=1)\cr
p(u_t=1|u_{t-2}=1,u_{t-1}=2)\cr
p(u_t=2|u_{t-2}=1,u_{t-1}=2)\cr
p(u_t=1|u_{t-2}=2,u_{t-1}=1)\cr
p(u_t=2|u_{t-2}=2,u_{t-1}=1)\cr
p(u_t=1|u_{t-2}=2,u_{t-1}=2)\cr
p(u_t=2|u_{t-2}=2,u_{t-1}=2)
\end{matrix}\right),\quad t=3,\ldots,T.
\]
Note that for $t=1$ this is a vector of initial probabilities and that
the number of elements of $\b p_t$ is $k^{d_{t,0}}$, where, in general,
$d_{t,j}=\min(t-1,h)+j+1$. Also let $\b f_t$ denote
the column vector with elements $f(y_t|u_t)$, $u_t=1,\ldots,k$, and let $\b q_{t,j}$ 
denote the column vector of the posterior probabilities 
$q(u_t|u_{\max(t-h,1)},\ldots,u_{t-1},u_{t+1},\ldots,u_{t+j},\b y)$
again arranged in lexicographical order. With $h=1$, $k=2$, and $j=1$, for instance, we have
\begin{equation}
\b q_{t,j}=\left(\begin{matrix}
q(u_t=1|u_{t-1}=1,u_{t+1}=1,\b y)\cr
q(u_t=1|u_{t-1}=1,u_{t+1}=2,\b y)\cr
q(u_t=2|u_{t-1}=1,u_{t+1}=1,\b y)\cr
q(u_t=2|u_{t-1}=1,u_{t+1}=2,\b y)\cr
q(u_t=1|u_{t-1}=2,u_{t+1}=1,\b y)\cr
q(u_t=1|u_{t-1}=2,u_{t+1}=2,\b y)\cr
q(u_t=2|u_{t-1}=2,u_{t+1}=1,\b y)\cr
q(u_t=2|u_{t-1}=2,u_{t+1}=2,\b y)
\end{matrix}\right),\quad t=2,\ldots,T.\label{eq:ex_vector_q}
\end{equation}
The dimension of this vector is $k^{d_{t,j}}$; note that for $j=0$ these
vectors contain the target posterior probabilities in (\ref{eq:def:p_t}). 

Finally, let $\b M_{a,b}$
be a marginalization matrix such that, given a column vector $\b v$ with elements indexed
by a sequence of $b$ variables assuming $k$ possible values
(as the above vectors), $\b M_{a,b}\b v$ provides the
corresponding vector in which the elements are summed with respect to the $a$-th of
these variables. This matrix may be simply constructed by following Kronecker product 
$\b M_{a,b}=\bigotimes_{l=1}^b \b M^*_l$, where
\[
\b M^*_l=\left\{
\begin{array}{ll}
\b 1_k\tr, & l=a,\\
\b I_k, & l\neq a,
\end{array}
\right.
\]
with $\b 1_k$ denoting a column vector of $k$ ones and $\b I_k$ denoting an identity
matrix of the same dimension. For instance, in the same context of the example that led
to the vector $\b q_t$ in (\ref{eq:ex_vector_q}), we have
\[
\b M_{2,3}\b q_{t,j}=\left(\begin{matrix}
q(u_t=1|u_{t-1}=1,u_{t+1}=1,\b y)+q(u_t=2|u_{t-1}=1,u_{t+1}=1,\b y)\cr
q(u_t=1|u_{t-1}=1,u_{t+1}=2,\b y)+q(u_t=2|u_{t-1}=1,u_{t+1}=2,\b y)\cr
q(u_t=1|u_{t-1}=2,u_{t+1}=1,\b y)+q(u_t=2|u_{t-1}=2,u_{t+1}=1,\b y)\cr
q(u_t=1|u_{t-1}=2,u_{t+1}=2,\b y)+q(u_t=2|u_{t-1}=2,u_{t+1}=2,\b y)
\end{matrix}\right).
\]

Using the above notation, for $t=T$ we directly obtain the target vector $\b q_{t,0}$
through the following operations, which directly derive from (\ref{eq:p_T}):
\begin{eqnarray*}
\b a_{T,0}&=&(\b 1_{k^{(d_{T,0}-1)}}\ot\b f_T)\times\b p_T,\\
\b q_{T,0} &=& \b a_T/(\b M_{d_{T,0},d_{T,0}}\tr\b M_{d_{T,0},d_{T,0}}\b a_{T,0}),
\end{eqnarray*}
where $\times$ and $/$ denote, respectively, elementwise product and division. Then,
for $t=1,\ldots,T-1$ (in reverse order), we first compute $\b q_{t,j}$ for
$j=\min(T-t,h)$ and then we recursively compute $\b q_{t,j}$ from
$j=\min(T-t,h)-1$ to $j=0$. In particular, for $j=\min(T-t,h)$ we  
compute the vector $\b a_{t,j}$ containing the elements at the numerator
of (\ref{eq:prob_cond2}) by the following recursion:
\begin{eqnarray*}
\b a_{t,0}&=&(\b 1_{k^{(d_{t,0}-1)}}\ot\b f_t)\times \b p_t,\\
\b a_{t,l}&=&(\b a_{t,l-1}\ot\b 1_k)\times(\b 1_{k^{(d_{t,l}-d_{t+l,0})}}
\ot\b p_{t+l}),\quad l=1,\ldots,j.
\end{eqnarray*}
Then, we have
\[
\b q_{t,j} = \b a_{t,j}/(\b M_{d_{t,0},d_{t,j}}\tr\b M_{d_{t,0},d_{t,j}}\b a_{t,j}).
\]
Finally, from (\ref{eq:theo}) we have the recursion
\begin{eqnarray*}
\b s_{t,{j+1}} &=& (\b 1_{k^{(d_{t,j+1}-d_{t+j+1,0})}}\ot\b q_{t+j+1,0})/\b q_{t,j+1},\\
\b q_{t,j} &=& \b 1_{k^{d_{t,j}}}/(\b M_{d_{t,j+1},d_{t,j+1}}\b s_{t,j+1}),
\end{eqnarray*}
to be applied for $j=\min(T-t,h)-1$ until $j=0$, when we obtain the target vector $\b q_{t,0}$.

In order to express the forward  recursion in (\ref{eq:forward_recursion2}) using the 
matrix notation, let $\b q_t^*$ denote the vector with elements 
$q(u_{\max(t-h,1)},\ldots,u_t|\b y)$ arranged in the usual lexicographical order.
Then, we have $\b q^*_1=\b q_{1,0}$, whereas for $t>1$, we have
\[
\b q^*_t=\left\{\begin{array}{ll}
\b q_{t,0}\times(\b q^*_{t-1}\ot\b 1_k), & t=2,\ldots,h+1,\\
\b q_{t,0}\times[(\b M_{1,h+1}\b q^*_{t-1})\ot\b 1_k], & t=h+2,\ldots,T.
\end{array}
\right.
\]

\bibliography{biblio}
\bibliographystyle{apalike}
\end{document}